\documentclass[12pt,twoside,reqno]{amsart}
\usepackage{nicefrac}
\usepackage[letterpaper, portrait, margin=1in]{geometry}
\usepackage[backend=biber,style=alphabetic,sorting=nyt]{biblatex}
\usepackage{amsmath}
\usepackage{amssymb}
\usepackage{fancyhdr}
\usepackage{mathrsfs}
\usepackage{indentfirst}
\usepackage{mathtools} 
\usepackage{microtype}
\usepackage{commath}
\usepackage[dvipsnames]{xcolor}
\usepackage{soul}
\usepackage{upgreek} 
\usepackage{graphicx} 
\graphicspath{ {images/} }
\usepackage[rightcaption]{sidecap}
\usepackage{wrapfig}
\usepackage{caption}
\usepackage{theoremref}
\usepackage{cancel}
\usepackage{amsthm}
\usepackage{xcolor}
\usepackage[final]{hyperref}
\usepackage{float}
\usepackage{comment}
\usepackage{soul}
\usepackage[normalem]{ulem}

\hypersetup{unicode= false, colorlinks=true, linkcolor=RoyalBlue,
	anchorcolor=RoyalBlue, citecolor=ForestGreen, filecolor=red, menucolor=RoyalBlue, urlcolor=RoyalBlue}

\newtheorem{theorem}{Theorem}[section]

\newtheorem{definition}[theorem]{Definition}
\newtheorem{example}[theorem]{Example}
\newtheorem{corollary}[theorem]{Corollary}
\newtheorem{remark}[theorem]{Remark}

\def\D{{\mathbb{D}}}
\def\C{{\mathbb{C}}}
\def\R{{\mathbb{R}}}

\def\HP{{\mathbb{C}_+}}

\newcommand{\stkout}[1]{\ifmmode\text{\sout{\ensuremath{#1}}}\else\sout{#1}\fi}

\def \nlhat {\overset{\ \curlywedge}}
\def \nlihat {\overset{\ \curlyvee}}

\title{Convergence from the discrete to the continuous non-linear Fourier transform}

\begin{document}

\setlength{\footskip}{14.0pt}.

\thispagestyle{empty}	

\author{Ashley R. Zhang}
	\address{Vanderbilt University\\ Department of Mathematics\\
		1326 Stevenson Center\\
		Nashville, TN  37240\\ USA }
	\email{ashley.zhang@vanderbilt.edu}
	
\begin{abstract}
    In this note, we study the convergence from the discrete to the continuous non-linear Fourier transform. Relations between spectral problems and questions in complex function theory provide a new approach to the study of scattering problems and the non-linear Fourier transform \cite{Scatter}. In particular, the non-linear Fourier transform can be viewed from the perspective of spectral problems for differential operators. Results in \cite{MP, PZ} can be seen as results for the non-linear Fourier transform. These results are similar to some convergence problems for the discrete non-linear Fourier transform considered in \cite{T} and \cite{TT}.
\end{abstract}
	
	\maketitle

\section{Introduction}

We study the convergence from the discrete to the continuous inverse non-linear Fourier transform in this note.

Over the past few decades, a wealth of evidence has indicated that scattering problems for differential operators can be seen as a non-linear version of the classical Fourier transform \cite{T, TT}.

Naturally, one might ask whether certain properties of the linear Fourier transform hold in the non-linear scattering setting. These problems have been studied over the past several decades and remain active today. Most of the developments and details of these problems can be found in \cite{Silva} and references therein. One of the more recent developments in the non-linear Fourier transform is pointwise convergence.

Pointwise convergence of the Fourier transform on $\R$, due to L. Carleson \cite{Carleson}, is one of the fundamental results in Fourier analysis. Poltoratski recently proved a non-linear analog of Carleson's theorem on the half-line $\R_+$ \cite{Scatter}. His approach is independent of existing results and is based on the study of resonances of Dirac systems using spectral theory for differential operators and complex function theory. Although it had been known that scattering and spectral problems are closely related to the non-linear Fourier transform, this connection was never made explicit until \cite{Scatter}.

The non-linear Fourier transform (NLFT) defined in \cite{TT} relies on a limit, similar to its linear counterpart. Moreover, the definition is only "explicit" for sequences. In general, there is no formula for the NLFT like there is for the linear Fourier transform. In section \ref{NLFT}, we propose an alternative definition for the NLFT, which depends on the spectral measure of the underlying Dirac (or canonical) system. With this definition, there is a natural way to construct a family of discrete measures whose NLFTs converge to the transform of a continuous function. Similar problems were considered in \cite{TT} without much detail. We are able to provide some of the details for this approach using our modified definition. 

The explicit formula for the discrete NLFT in \cite{T} involves a product of matrix exponentials, which are difficult to compute. Even for a finite sequence, it is challenging to obtain an explicit example using the formula. However, as far as we know, this sequential formula is the only example of the NLFT in the literature. Our modified definition uses the spectral measure of a Dirac system and avoids explicitly taking the limit. With this definition, we are able to connect some sequences with their NLFTs and obtain examples of the discrete NLFT.

The main object in \cite{TT} is the discrete NLFT. It was mentioned there (without proof) that the NLFTs of discrete sequences can be used to approximate the NLFTs of functions. For a Dirac system with evenly spaced discrete measures as "potentials", the smaller the distance between adjacent pointmasses is, the larger the period of its spectral measure becomes. Using the periodization approach introduced in \cite{PZ}, we are able to construct sequences (discrete measures) whose NLFTs converge to the NLFT of a function.

The paper is organized as follows. We introduce real Dirac systems and canonical systems in section \ref{RDS-CS}, the basics of inner function theory and Krein-de Branges theory in section \ref{kdB}, and spectral measures in section \ref{SpectralMeasure}. The definitions for the non-linear Fourier transform of locally integrable functions and discrete measures, as well as some of their properties, are contained in section \ref{NLFT}. Our main result, using the NLFTs of discrete measures to approximate the NLFT of continuous functions, is presented in section \ref{results}. Finally, in section \ref{examples}, we provide several explicit examples of the NLFT, which are less frequently documented in the existing literature.

\section*{Acknowledgement}
This work is an application of joint work \cite{PZ} with Alexei Poltoratski, from whom I learned spectral problems and the non-linear Fourier transform. I am also indebted for all the help and support he provided during my graduate studies.

\section{Real Dirac Systems and Canonical Systems}\label{RDS-CS}

Real Dirac systems on the right half-line $\mathbb{R}_+$ are systems of the form \begin{equation}\label{RD}
    \Omega \Dot{X}(t) = z X(t) - Q(t) X(t), \ \ t \in (0, \infty),
\end{equation}
where $z\in \mathbb{C}$ is a spectral parameter,
$\Omega = \begin{pmatrix} 0 & 1 \\ -1 & 0 \end{pmatrix}$ is the symplectic matrix, and $Q(t) = \begin{pmatrix} 0 & f(t) \\ f(t) & 0 \end{pmatrix}$ is the potential matrix, where $f(t)$ is a real-valued locally summable function, and $X(t) = X(t,z) =  \begin{pmatrix} u(t,z) \\ v(t,z) \end{pmatrix}$ is the unknown vector function.

Real Dirac systems can be rewritten as canonical systems, which form a broad class of second-order differential equation systems. Canonical systems are $2 \times 2$ differential equation systems of the form \begin{equation} \label{CS}
    \Omega \Dot{X}(t) = z H(t) X(t), \ \ t \in (t_-, t_+),\ \ -\infty < t_- < t_+ \leqslant \infty,
\end{equation}
where $z \in \C$ is a spectral parameter,
$\Omega = \begin{pmatrix} 0 & 1 \\ -1 & 0 \end{pmatrix}$ is the symplectic matrix as before,
$H(t) = \begin{pmatrix} h_{11}(t) & h_{12}(t) \\ h_{21}(t) & h_{22}(t)\end{pmatrix}$ is a given matrix-valued function called the Hamiltonian of the system, and $X(t) = X(t,z) = \begin{pmatrix} u(t,z) \\ v(t,z) \end{pmatrix}$ is the unknown vector function. Here $H(t) \in \mathbb{R}^{2 \times 2}$, $H(t) \in L^1_{loc}(\mathbb{R})$, and $H(t) \geqslant 0$ almost everywhere.

For the reader's convenience, we provide details on rewriting real Dirac systems as canonical systems. Details of this process can be found for instance in \cite{MPS}.

Let $V$ be the real $2 \times 2$ matrix satisfying \begin{equation*}
\begin{cases}
\Omega \Dot{V}(t) = - Q(t) V(t) \\
V(0) = I
\end{cases}.  
\end{equation*}

The solution to the system above is \begin{equation*}
    V(t) = \begin{pmatrix} \exp(\int_0^t f(x) dx) & 0 \\ 0 &  \exp(-\int_0^t f(x) dx)\end{pmatrix}.
\end{equation*}

We claim that this matrix $V(t)$ satisfies $V^\ast(t) \Omega V(t) = \Omega$ for all $t$. First note that $V^\ast(0) \Omega V(0) = \Omega$. It remains to show that $V^\ast(t) \Omega V(t)$ stays constant. We differentiate $V^\ast(t) \Omega V(t)$ and obtain \begin{equation*}
\begin{split}
    \Dot{V^\ast}(t) \Omega V(t) + V^\ast(t) \Omega \Dot{V}(t) &= -\left(\Omega \Dot{V}(t)\right)^\ast V(t) + V^\ast(t) \left(-Q(t) V(t)\right)\\
    & = \left(Q(t) V(t) \right)^\ast V(t) - V^\ast(t) Q(t) V(t) = 0.
\end{split}
\end{equation*}

Finally, with a change of variable $X(t) = V(t) Y(t)$, where $X(t)$ is the vector solution in \eqref{RD}, we get \begin{equation*}
\begin{split}
    &\Omega \Dot{\left( V(t) Y(t) \right)} = z V(t) Y(t) - Q(t) V(t) Y(t), \\
    &\Omega \Dot{V}(t) Y(t) + \Omega V(t) \Dot{Y}(t) = z V(t) Y(t) - Q(t) V(t) Y(t), \\
    & \cancel{\left( - Q(t) V(t) \right) Y(t)} + \Omega V(t) \Dot{Y}(t) = z V(t) Y(t) \cancel{- Q(t) V(t) Y(t)}, \\
    & V^\ast  \Omega V(t) \Dot{Y}(t) = z  V^\ast(t) V(t) Y(t),
\end{split}
\end{equation*}
which is a canonical system for $Y(t)$ with Hamiltonian \begin{equation}\label{RDHamiltonian}
    H^D = V^\ast(t) V(t) = \begin{pmatrix} \exp \left( 2 \int_0^t f(x)dx \right) & 0 \\ 0 & \exp \left( -2\int_0^t f(x)dx \right) \end{pmatrix}.
\end{equation}

If distributions are allowed in $f$, all canonical systems with diagonal determinant-normalized Hamiltonians can be brought to real Dirac form as well. In other words, we broaden the class of $f$ we allow in \eqref{RD}: any real $f$ satisfying $\exp(2 \int_0^t f(x) dx) \in L^1_{loc}(\mathbb{R}_+)$ can be the potential of a Dirac system. In particular, we allow discrete measures to be potentials of Dirac systems.  Krein-de Branges theory for canonical systems applies to real Dirac systems, see \cite{dB, Dym, Remling, Rom}. The next section contains a brief summary of the theory that we will need.

\section{Inner functions and Krein-de Branges theory}\label{kdB}

\subsection{Inner functions and Clark measures}
This brief discussion on inner functions and Clark measures is primarily based on \cite{CBMS} and \cite{PS}.

We use $H^\infty(\HP)$ to denote the space of bounded analytic functions on the upper half-plane $\HP$. An \textit{inner function} is a function $\Theta$ in $H^\infty(\HP)$ satisfying $\abs{\Theta(x)} = 1$ almost everywhere on $\R$. 

Analytic functions with non-negative imaginary parts are called \textit{Herglotz functions}. We say a measure $\mu$ is Poisson-finite if $$\int \frac{d\abs{\mu(t)}}{1 + t^2} < \infty.$$

\begin{theorem} [Herglotz representation theorem for $\HP$]
Let $m: \HP \to \HP$ be analytic. Then there exists a positive Poisson-finite measure $\mu$ on $\R$ and a constant $A \geqslant 0$ such that \begin{equation}\label{HerglotzRep}
    m(z) = \Re{m(i)} + Az + \frac{1}{\pi} \int \left( \frac{1}{t - z} - \frac{t}{1 + t^2} \right) d\mu(t).
\end{equation}
\end{theorem}

The Cayley transform given by $z \mapsto i\frac{1 + z}{1 - z}$ maps the unit disk $\mathbb{D}$ conformally to the upper half-plane $\mathbb{C}_+$. Thus, for every inner function $\Theta(z)$, there is a corresponding Herglotz function, \begin{equation*}
    m(z) = i \frac{1 + \Theta(z)}{1 - \Theta(z)},
\end{equation*}
obtained by composing the inner function with the Cayley transform.

Given a Herglotz function $m(z)$, the corresponding $\Theta$ is inner if and only if the measure $\mu$ in \eqref{HerglotzRep} is singular with respect to the Lebesgue measure. Since we will use this to define the Clark measure, we write it down explicitly \begin{equation}\label{ClarkDef}
    \Re \frac{1 + \Theta(z)}{1 - \Theta(z)} = A y + \frac{1}{\pi} \int \frac{y d\mu(t)}{(x - t)^2 + y^2},
\end{equation}
where $A$ can be interpreted as a point mass at infinity, and the measure $\mu$ is singular with respect to the Lebesgue measure, Poisson-finite, and supported on the set of $\R$ where the non-tangential limits of $\Theta(z)$ are $1$. The measure $\sigma = \mu + A \delta_\infty$ on $\hat{\mathbb{R}} = \mathbb{R} \cup \{ \infty \}$ is called the Clark measure of $\Theta$.

For a given inner function $\Theta(z)$, we can define a family of Clark measures. The measure we just obtained is usually denoted by $\sigma_1$. Let $\alpha \in \mathbb{T}$ be a unimodular number, then $\Tilde{\Theta}(z) = \bar{\alpha} \Theta(z)$ is still an inner function and we can obtain a different Clark measure $\Tilde{\sigma}$ from \eqref{ClarkDef}. This Clark measure is usually denoted by $\sigma_\alpha$. We call the family of measures $\{\sigma_\alpha\}_{\abs{\alpha} = 1}$ the family of Clark measures of an inner function $\Theta$.

Meromorphic inner functions (MIF) are inner functions in the upper half-plane $\mathbb{C}_+$ that can be extended meromorphically to the whole complex plane. In this case, we can write down the Clark measure $\sigma_\alpha$ explicitly \begin{equation*}
    \sigma_{\alpha} = \frac{2\pi}{\abs{\Theta'(x)}} \delta_x, \quad \text{where~} x \in \{\Theta = \alpha\}.
\end{equation*}

\subsection{Model spaces}
We use $H^2(\mathbb{C}_+)$ to denote the Hardy space on the upper half-plane $\mathbb{C}_+$: \begin{equation*}
    H^2(\mathbb{C}_+) = \{ f \in \mathcal{H}(\mathbb{C_+})~|~ \sup_{y > 0} \int_{\mathbb{R}} \abs{f(x + iy)}^2 dx < \infty\},
\end{equation*}
where $\mathcal{H}(\mathbb{C}_+)$ denotes the set of all analytic functions in $\mathbb{C_+}$.

For each inner function $\Theta$, we can consider the associated model space \begin{equation*}
    K_\Theta = H^2(\HP) \ominus \Theta H^2(\HP),
\end{equation*}
where $\Theta H^2(\HP) = \left\{ \Theta(z) f(z) | f(z) \in H^2(\HP) \right\}$. Here $\ominus$ is the orthogonal difference between the two subspaces.

$K_\Theta$ is a reproducing kernel Hilbert space. The reproducing kernel is \begin{equation*}
    k_\Theta(z, w) = \frac{1}{2 \pi i} \cdot \frac{1 - \Theta(z) \overline{\Theta(w)}}{\overline{w} - z},
\end{equation*}
and for every $f \in K_\Theta$, $f(z) = \langle f, k_\Theta \rangle_w$.

Recall that the Clark measure of the inner function $\Theta$ is denoted by $\sigma$. Functions in $K_\Theta$ have non-tangential boundary values $\sigma-$almost everywhere, and can be recovered from the boundary values using the following formula \begin{equation}\label{ModelBdry}
    f(z) = \frac{1}{2 \pi i} \left( 1 - \Theta(z) \right) \int f(t) \overline{\left( 1 - \Theta(t) \right)} dt + \frac{1 - \Theta(z)}{2 \pi i} \int \frac{f(t)}{t - z} d\sigma(t).
\end{equation}

\subsection{Entire functions and de Branges spaces}
An entire function $E(z)$ is called an \textit{Hermite-Biehler function} if $|E(z)| > |E(\overline{z})|$ for all $z \in \mathbb{C}_+$.
Given an entire function $E(z)$ of Hermite-Biehler class, define the \textit{de Branges} space $B(E)$ based on $E$ as \begin{equation*}
    B(E) = \left\{F \text{~entire}\ |\ F/E, F^\#/E \in H^2(\mathbb{C_+}) \right\},
\end{equation*}
where $F^\#(z) = \overline{F(\overline{z})}$.  The space $B(E)$ becomes a Hilbert space when endowed with the scalar product inherited from $H^2(\C_+)$: \begin{equation*}
    [F, G] = \int_{-\infty}^\infty F(t) \overline{G(t)} \frac{dt}{|E(t)^2|}.
\end{equation*}We call an entire function real if it is real-valued on $\R$. Associated with $E$ are two real entire functions  $A = \frac{1}{2}(E + E^\#)$ and $C = \frac{i}{2}(E - E^\#)$ such that $E=A-iC$.
De Branges spaces are also reproducing kernel Hilbert spaces. The reproducing kernels 
\begin{equation*}
    K_\lambda(z) = \frac{\overline{E(\lambda)} E(z) - \overline{E^\#(\lambda)} E^\#(z)}{2\pi i(\overline{\lambda} - z)} = \frac{\overline{C(\lambda)} A(z) - \overline{A(\lambda)}C(z)}{\pi(\overline{\lambda} - z)}
\end{equation*}
are functions from $B(E)$ such that $[F, K_{\lambda}] = F(\lambda)$ for all $F \in B(E)$, $\lambda \in \mathbb{C}$.

The most elementary example of de Branges spaces is the Paley-Wiener spaces. We denote by $PW_t$ the standard Paley-Wiener space: the subspace of $L^2(\R)$ consisting of entire functions of exponential type at most $t$. $PW_t=B(E)$ with $E(z)=e^{-itz}$.
The reproducing kernels of $PW_t$ are the sinc functions,
$$\frac{\sin \left( t (z-\lambda) \right) }{\pi( z - \lambda)}.$$

De Branges spaces have an alternative axiomatic definition, which is useful in many applications.

\begin{theorem}[\cite{dB}] \label{dBAxioms}
Suppose that $H$ is a Hilbert space of entire functions that satisfies the following conditions: \begin{itemize}
\item $F \in H$, $F(\lambda) = 0$ $\Rightarrow$ $F(z) \frac{z - \overline{\lambda}}{z - \lambda} \in H$ with the same norm,
\item $\forall \lambda \notin \mathbb{R}$, the point evaluation is a bounded linear functional on $H$,
\item The map $F \to F^\#$ is an isometry.
\end{itemize}
Then $H = B(E)$ for some entire function $E$ of Hermite-Biehler class.
\end{theorem}

Every Hermite-Biehler entire function $E(z)$ gives rise to a meromorphic inner function $\Theta(z) = E^\#(z)/E(z)$, and therefore a corresponding model space $K_\Theta$. There exists an isomorphism between $K_\Theta$ and $B(E)$ given by $F \mapsto EF$. 

Conversely, for a meromorphic inner function $\Theta$, there exists some Hermite-Biehler entire function $E$ such that $\Theta(z) = E^\#(z)/E(z)$. As mentioned before, the Clark measure of a meromorphic inner function $\Theta$ is a discrete measure $\sigma_{\alpha} = \frac{2\pi}{\abs{\Theta'(x)}} \delta_x$ where $x \in \{\Theta = \alpha\}$, and $\alpha \in \mathbb{T}$. We will call the measure $\mu_\alpha = \abs{E}^2 \sigma_{\alpha}$ the \textit{representing measure} of the de Branges space $B(E)$. We can write down the measure $\mu_\alpha$ explicitly as well: \begin{equation*}
    \mu_\alpha = \frac{\pi}{- A'C + A C'} \delta_x,
\end{equation*} where $x \in \{ \frac{E^\#}{E} = \alpha \}$. We will pay special attention to $\mu_1$ and $\mu_{-1}$. Since $\frac{E^{\#}}{E} = 1$ when $C(z) = 0$ and $\frac{E^{\#}}{E} = -1$ when $A(z) = 0$, we have \begin{equation*}
    \mu_{1} = \frac{\pi}{A C'} \delta_x, \quad x \in Z(C),
\end{equation*}
and \begin{equation*}
    \mu_{-1} = \frac{\pi}{A' C} \delta_x, \quad x \in Z(A).
\end{equation*}
Here $Z(C)$ and $Z(A)$ stand for the zero sets of the functions $C$ and $A$, respectively.

\section{Spectral Measures}\label{SpectralMeasure}
\subsection{Spectral measures}
A solution to \eqref{RD} is a $C^2-$function 
$$X(t, z) = \begin{pmatrix} u(t, z) \\ v(t, z) \end{pmatrix}$$
on $(0, \infty)$ satisfying the equation system. An initial value problem (IVP) for \eqref{RD} is given by an initial condition $X(0) = c$, $c \in \mathbb{R}^2$. 
It follows from general existence and uniqueness results for ODEs that every IVP has a unique solution, see for instance chapter 1 of \cite{Remling}.

Let $X(t, z) = \begin{pmatrix} u(t, z) \\v(t, z) \end{pmatrix}$ be the unique solution for \eqref{RD} on $(0, \infty)$
satisfying a self-adjoint initial condition. Krein first observed that the function $E(t, z) = u(t, z) - i v(t, z)$ is an Hermite-Biehler entire function for every fixed $t$. 

Every real Dirac system delivers a family of nested de Branges spaces $B \left( E(t, z) \right)$: if  $0 < t_1 < t_2 \leqslant \infty$, then $B \left( E(t_1, z) \right) \subseteq B \left( E(t_2, z) \right)$, and all such inclusions are isometric because real Dirac systems do not have jump intervals (subintervals of $(0, \infty)$ on which the Hamiltonian in the canonical system is a constant matrix of rank $1$).

A positive measure $\mu$ on $\R$  is called the spectral measure of \eqref{RD} corresponding to the Neumann condition $\begin{pmatrix}    1 \\ 0\end{pmatrix}$ at $t = 0$ if all de Branges spaces $B\left( E(t, z) \right)$, $t \in (0, \infty)$ are isometrically embedded into $L^2(\mu)$.  
A spectral measure always exists, and it is unique if and only if \begin{equation*}
    \int_{0}^{\infty} \text{trace~} H(t) dt = \infty.
\end{equation*}
The case when the integral above is infinite (and the spectral measure is unique) is called the limit point case, otherwise it is called the limit circle case. 

Under our assumptions for real Dirac systems, any spectral measure $\mu$ is Poisson-finite; see for instance \cite{dB, Remling}. 

We will pay special attention to systems whose spectral measures are sampling for all Paley-Wiener spaces $PW_t$.

\begin{definition}
A positive Poisson-finite measure $\mu$ is sampling for the Paley-Wiener space $PW_t$ if there exist constants $0 < c < C$ such that for all $f \in PW_t$ \begin{equation*}
    c \|f\|_{PW_t} \leqslant \|f\|_{L^2(\mu)} \leqslant C \|f\|_{PW_t}.
\end{equation*}

A positive  measure $\mu$ is a Paley-Wiener measure ($\mu \in PW$) if it is sampling for $PW_t$ for all $t > 0$.
\end{definition}

It is not difficult to show that any PW-measure is Poisson-finite. For this and other basic properties of PW-measures, see \cite{MP} and references therein.A canonical system on the half-line $\mathbb{R}_+$ is called a PW-system if the corresponding de Branges spaces $B\left( E(t, z) \right)$ are  equal to  $PW_t$ as sets for all $t \in (0, \infty)$. In this case any spectral measure $\mu$ of the system is a PW-measure. It follows that PW-systems can be equivalently defined as those systems whose spectral measures belong to the PW-class. 

Real Dirac systems with locally summable potentials, as a particular subclass of canonical systems with diagonal determinant-normalized Hamiltonians, always give rise to PW-spectral measures. By \eqref{RDHamiltonian}, when an evenly spaced discrete measure serves as the potential of a real Dirac system, the Hamiltonian of the corresponding canonical system becomes a step function of uniform step size, and therefore the corresponding spectral measures are periodic, see theorem \ref{2piThm} and corollary \ref{2Talgorithm} below. 

The class of PW-measures admits the following elementary description. 

\begin{definition}
Let $\mu$ be a positive measure on $\R$. We call an interval $I \subseteq \mathbb{R}$ a $(\mu, \delta)$-interval if \begin{equation*}
    \mu(I) > \delta, \quad \text{and} \quad \abs{I} > \delta,
\end{equation*}
where $\abs{I}$ stands for the length of the interval $I$.
\end{definition}

\begin{theorem}\label{PWcondition}\cite{MP}
A positive Poisson-finite measure $\mu$ is a Paley-Wiener measure if and only if \begin{itemize}
    \item $\sup_{x \in \mathbb{R}} (x, x + 1) < \infty$,
    \item For any $d > 0$, there exists $\delta > 0$ such that for all sufficiently large intervals $I$, there exist at least $d \abs{I}$ disjoint $(\mu, \delta)$-intervals intersecting $I$.
\end{itemize}

\end{theorem}

We say that a measure $\mu$ has locally infinite support if $\text{supp~}\mu\cap [-C,C]$ is an infinite
set for some $C>0$, or equivalently if $\text{supp\ }\mu$ has a finite accumulation point.
For a periodic measure, one can easily deduce the following.

\begin{corollary}\label{PWper}\cite{MP}
A positive locally finite periodic measure is a Paley-Wiener measure if and only if it has locally infinite support.
\end{corollary}

\subsection{Hilbert transform}
This section is primarily based on section 19 of \cite{MPS}, and this is one of the key ingredients for our alternative definition for the non-linear Fourier transform in section \ref{NLFT}.

Consider the matrix initial value problem \begin{equation}\label{RDMatrixInf}
    \Omega \Dot{M}(t) = z M(t) - Q(t) M(t), \quad M(0) = I,
\end{equation}
on the half-line $\mathbb{R}_+$ and an interval $(0, b)$.

The matrix $M(t, z)$ is entire with respect to $z$ for every fixed $t$. The first column is the solution to the IVP \eqref{RD} satisfying the Neumann boundary condition at $t = 0$, and the second column is the solution satisfying the Dirichlet boundary condition at $t = 0$.

\begin{definition}
A matrix $\begin{pmatrix} A(z) & B(z) \\ C(z) & D(z) \end{pmatrix}$ is called a Nevanlinna matrix if all elements are entire functions and the following conditions are satisfied: \begin{itemize}
    \item this identity holds: \begin{equation*}
        A(z)D(z) - B(z)C(z) = 1;
    \end{equation*}
    \item for any fixed real $s$, the function \begin{equation*}
        w(z) = \frac{A(z) s - B(z)}{C(z)s - D(z)}
    \end{equation*}
    satisfies the inequality \begin{equation*}
        \frac{\Im{w(z)}}{\Im{z}} > 0.
    \end{equation*}
\end{itemize}
\end{definition}

For every fixed $t$, the matrix solutions to the initial value problems \ref{RDMatrixInf} on the half-line $\mathbb{R}_+$ and an interval $(0, b)$ are Nevanlinna matrices. In the rest of this section, we focus on the case on an interval $(0, b)$. For convenience, we will use $A(z)$ instead of $A(b, z)$, and similarly for $B$, $C$ and $D$. 

We consider the following objects: \begin{enumerate}
\item $E(z) = A(z) - iC(z)$ is the first Hermite-Biehler entire function, and $\mu_A$ is the representing measure of $B(E)$ supported on $Z(A)$, $\mu_C$ is its representing measure supported on $Z(C)$. Recall from section \ref{kdB} that $Z(A)$ and $Z(C)$ stand for the zero sets of functions $A(z)$ and $C(z)$ on $\R$. 

\item $F(z) = A(z) + iB(z)$ is our second Hermite-Biehler entire function. We consider the associated inner function \begin{equation*}
    \Theta(z) = \frac{A(z) - iB(z)}{A(z) + iB(z)},
\end{equation*}
and denote its Clark measures by $\sigma_\alpha$. Then $\sigma_{-1}$ is supported on $Z(A)$, and $\sigma_{1}$ is supported on $Z(B)$.
\end{enumerate}

Let $\mu^F$ be the representing measure of $B(F)$ supported on $Z(A)$. Recall that \begin{equation*}
    \mu_A = - \frac{\pi }{A' C}\delta_x, \quad \quad \mu^F = \frac{\pi }{A' B}\delta_x, \quad \quad x \in Z(A).
\end{equation*}

On $Z(A)$, since the matrix $M$ is normalized by its determinant, we have $B(z) C(z) = -1$, which implies $B(z) = - 1/C(z)$. Also note that on $Z(A)$, $\abs{E}^2 = \abs{C}^2$ and $\abs{F}^2 = \abs{B}^2$. Therefore, \begin{equation*}
    \mu_A = - \frac{\pi\delta_x}{A'} \frac{1}{C} = \frac{\pi \delta_x}{A'} B = \abs{B}^2 \frac{\pi \delta_x}{A' B} = \abs{F}^2 \mu_F = \sigma_{-1}.
\end{equation*}

Let $m$ be a Herglotz function, and $\Theta$ be the associated inner function obtained by composing the Herglotz function with the Cayley transform. The Schwarz transform, \begin{equation*}
    \mathcal{S}_\mu(z) = \frac{1}{\pi i} \int_\mathbb{R} \frac{1}{t - z} - \frac{t}{1 + t^2} d\mu(t),
\end{equation*} 
recovers a Herglotz function $m(z)$ satisfying $\Im m(i) = 0$ from the Clark measure of $\Theta$. Therefore, we can recover the inner function from its Clark measure up to a M\"obius transform. This means that for every finite $t \in (0, b)$, we have \begin{equation}\label{ratio}
    \frac{B(z)}{A(z)} = -\frac{i}{\mathcal{S}_{\sigma_{-1}}} = - \frac{i}{\mathcal{S}_{\mu_A}}.
\end{equation}

Similarly, we can look at the Hermite-Biehler functions $\Tilde{E}(z) = B(z) - i D(z)$ and $\Tilde{F}(z) = C(z) + i D(z)$, perform a similar calculation and obtain \begin{equation}\label{ratio'}
    \frac{D(z)}{C(z)} = -\frac{i}{\mathcal{S}_{\Tilde{\sigma}_{-1}}} = - \frac{i}{\mathcal{S}_{\mu_C}},
\end{equation}
where $\Tilde{\sigma}_{-1}$ is the Clark measure of the inner function $$ \Tilde{\Theta}(z) = \frac{C(z) - i D(z)}{C(z) + i D(z)}.$$

\section{The Non-linear Fourier Transform}\label{NLFT}

For the reader's convenience, we provide details on the connection between Hermite-Biehler entire functions and the non-linear Fourier transform discussed in \cite{Scatter}. Denote by $E(t, z)$ and $\tilde{E}(t, z)$ the Hermite-Biehler functions corresponding to the Neumann $X(0, z) = \begin{pmatrix}1 \\ 0\end{pmatrix}$ and Dirichlet $X(0, z) = \begin{pmatrix}0 \\ 1\end{pmatrix}$ boundary conditions, respectively. Then \eqref{RD} for $E$ and $\tilde{E}$ becomes \begin{equation*}
    \frac{\partial}{\partial t} E(t, z) = \left( - z v(t, z) + f(t) u(t, z) \right) - i \left( z u(t, z) - f(t) v(t, z) \right),
\end{equation*}
which yields \begin{equation}\label{RDHB}
    \frac{\partial}{\partial t} E(t, z) = -iz E(t, z) + f(t) E^\#(t, z),
\end{equation}
with initial conditions $E(0, z) = 1 - 0 i = 1$ for $E$, and $\tilde{E}(0, z) = 0 - i = -i$ for $\tilde{E}$.

For each $t \geqslant 0$, define the following entire functions \begin{equation}\label{pair}
\begin{split}
    a(t, z) &= \frac{\exp(i t z)}{2} \left( E(t, z) + i \tilde{E}(t, z) \right),\\
    b(t, z) &= \frac{\exp(i t z)}{2} \left( E(t, z) - i \tilde{E}(t, z) \right).
\end{split}    
\end{equation}
As in \cite{Scatter}, this notation is slightly different from the one in \cite{T} and \cite{TT}, where their $a$ is $a^\#$ here.

The limit of $\left(a(t, z), b(t, z)\right)$ as $t \to \infty$ can be viewed as a non-linear analog of the Fourier transform if the limit exists in some sense. We denote this limit by $\left( a(\infty,z), b(\infty,z)\right)$. This definition was used for instance in \cite{T}. It is known that when $f \in L^2$, the pair $\left( a(t, z), b(t, z) \right)$ converges normally on the half-plane $\HP$ (see \cite{Denisov}), and the function $b(t, z) / a(t, z)$ converges pointwise on the half-line $\R_+$ (see \cite{Scatter}).

Using \eqref{RDHB}, one can show that the matrix \begin{equation*}
    G(t, z) = \begin{pmatrix} a^\#(t, z) & b^\#(t, z) \\ b(t, z) &  a(t, z) \end{pmatrix}
\end{equation*}
satisfies the differential equation \begin{equation}\label{DEqDef}
    \frac{\partial }{\partial t} G(t, z) = \begin{pmatrix}0 & e^{-2izt} f(t) \\ e^{2izt} f(t) & 0 \end{pmatrix} G(t, z),
\end{equation}
with initial condition $G(0, z) = I$. 

If we study this \eqref{DEqDef} on its own, it is not hard to show that the solution to this IVP is of the form \begin{equation*}
    G(t, z) =  \begin{pmatrix} a^\#(t, z) & b^\#(t, z) \\ b(t, z) &  a(t, z) \end{pmatrix}.
\end{equation*}
The limit of the solution $\left( a(t, z), b(t, z) \right)$ as $t \to \infty$ is analogous to the linear Fourier transform if the limit exists in some sense. This is the definition used, for instance, in \cite{MTT}.

As in \cite{Scatter}, we will pay special attention to the functions \begin{equation*}
    \nlhat{f}(z) = \frac{b(\infty, z)}{a(\infty, z)},
\end{equation*} and  \begin{equation*}
    \nlhat{f}_t(z) = \frac{b(t, z)}{a(t, z)}.
\end{equation*}
The $\nlhat{f}_t$'s contain enough information to recover $(a(t, z), b(t,z))$: $a(t, z)$ can be uniquely recovered since $a(t, z)$ is outer for all $t$, positive at $0$, and the absolute value of $a(t, z)$ satisfies \begin{equation*}
    |a(t, z)| = \frac{1}{\sqrt{1 - |\nlhat{f}_t|^2}} \quad \text{ on } \mathbb{R},
\end{equation*}
and $b(t, z) = a(t, z) \nlhat{f}_t$. Under the assumption $f \in L^2$, we have $\nlhat{f}_t \to \nlhat{f}$ almost everywhere \cite{Scatter}.

This function $\nlhat{f}$ also has a nice connection with the spectral measure of the underlying Dirac system: by \eqref{pair}, for every finite $t$, we have \begin{equation*}
    \frac{b(t, z)}{a(t, z)} = \frac{E(t, z) - i \tilde{E}(t, z)}{E(t, z) + i \tilde{E}(t, z)} = \frac{1 - i\frac{\tilde{E}(t, z)}{E(t, z)}}{1 + i\frac{\tilde{E}(t, z)}{E(t, z)}}.
\end{equation*}

Recall from the previous section that $\mu_A$ and $\mu_C$ are the representing measures of $B(E)$, supported on $Z(A)$ and $Z(C)$, respectively. Since all the de Branges spaces are isometrically embedded into $L^2(\mu)$, both $\mu_A$ and $\mu_C$ tend to $\mu$ as $t \to \infty$ $\ast-$weakly on $\hat{\mathbb{R}} = \mathbb{R} \cup \{\infty\}$. Together with \eqref{ratio}, we get \begin{equation*}
    \frac{\tilde{u}(t, z)}{u(t, z)}, \frac{\tilde{v}(t, z)}{v(t, z)} \to -i\frac{1}{\mathcal{S}_\mu(z)}, \quad \text{as} \; t \to \infty,
\end{equation*}
where $\mathcal{S}_\mu(z)$ is the Schwarz transform.

Then we have \begin{equation*}
    \frac{b(t, z)}{a(t, z)} = \frac{1 - i\frac{\tilde{E}(t, z)}{E(t, z)}}{1 + i\frac{\tilde{E}(t, z)}{E(t, z)}} \to \frac{\mathcal{S}_\mu - 1}{\mathcal{S}_\mu + 1}.
\end{equation*}

Based on this, we propose the following definition for the non-linear Fourier transform.

\begin{definition}\label{NLFTDef}
Let $f$ be a real-valued locally summable function on the half-line $\mathbb{R}_+$, and $\mu$ be the spectral measure of the real Dirac system with potential $Q = \begin{pmatrix} 0 & f \\ f & 0 \end{pmatrix}$. The non-linear Fourier transform (NLFT) of $f$ is \begin{equation*}
    \nlhat{f}(z) = \frac{\mathcal{S}_\mu - 1}{\mathcal{S}_\mu + 1},
\end{equation*}
where $\mathcal{S}_\mu(z)$ is the Schwarz transform.
\end{definition}

One advantage of this definition is that it avoids explicitly taking the limit of $\left( a(t, z), b(t, z) \right)$, which is only guaranteed to exist in the sense of normal convergence on the half-plane $\mathbb{C}_+$ but may not in the strong sense on $\mathbb{R}$. The definition above agrees with $\lim_{t \to \infty} \left( a(t, z), b(t, z) \right)$ if such a limit exists in some sense. This definition also applies to the discrete case studied by Tao and Thiele in \cite{TT} if discrete evenly spaced distributions are allowed to be the potential $f$. We will discuss the discrete case further in the second half of this section.

This transform $f \to \nlhat{f}$ has several properties similar to the linear Fourier transform. The following properties are summarized in \cite{Silva}:
\begin{itemize}
\item Modulation/ frequency shifting: if $h(t) = e^{2 i s t} {f}(t)$, then $\nlhat{h}(z) = \nlhat{f}(z - s)$ \footnote{Since we only consider real potentials, this property does not apply to us. This is generally true when complex potentials are considered.};

\item Translation / time shifting: if $h(t) = f(t - t_0)$, then $\nlhat{h}(z) = e^{2 i t_0 z}\nlhat{f}(z)$;

\item Time scaling: for $a > 0$ and $h(t) = a f(at)$, then $\nlhat{h}(z) = \nlhat{f}(\frac{z}{a})$;

\item Linearization at the origin is the Fourier transform \footnote{Here the Fourier transform is given by $\hat{f}(\xi) = \int_{-\infty}^\infty \exp(2 i \xi x) f(x) dx$.}: \begin{equation*}
    \nlhat{\epsilon f} (z) = \epsilon \hat{f}(z) + \mathcal{O}(\epsilon^2), \quad \text{~as~} \epsilon \to 0;
\end{equation*}

\item Non-linear Riemann-Lebesgue: \begin{equation*}
    |\nlhat{f}(z)| \leqslant \tanh \left( \norm{f}_{L^1(\mathbb{R})} \right);
\end{equation*}

\item Non-linear Parseval: \begin{equation*}
    \|f\|_{L^2(\mathbb{R})}^2 = \frac{1}{2} \norm{ \log \left(1 - |\nlhat{f}|^2 \right)}_{L^1(\mathbb{R})}.
\end{equation*}
\end{itemize}

Tao and Thiele primarily studied the non-linear Fourier transform of discrete sequences in \cite{TT}. The sequential case can be studied via real Dirac systems if we choose discrete measures to be the potential. Let $f$ be an evenly spaced discrete measure
\begin{equation*}
    f(t) = \sum_{n = 1}^\infty c_n \delta_{n}(t), \quad \quad c_n \in \mathbb{R}.
\end{equation*}
Define $a(N, z)$ and $b(N, z)$ iteratively as follows:  \begin{equation*}
\begin{split}
    a(0, z) &= 1, \quad b(0, z) = 0, \\
    \begin{pmatrix} a^\#(N, z) & b^\#(N, z) \\ b(N, z) &  a(N, z) \end{pmatrix} &=   \exp \begin{pmatrix} 0 & c_N \exp( - 2 i z N) \\ c_N \exp(2 i z N) & 0 \end{pmatrix} \begin{pmatrix} a^\#(N - 1, z) & b^\#(N - 1, z) \\ b(N - 1, z) &  a(N - 1, z) \end{pmatrix}.    
\end{split}
\end{equation*}

Note that this is different from the discrete transform studied in \cite{TT}, since they worked with $c_n \in \mathbb{T}$. This recursion above can be viewed as a discrete version of \eqref{DEqDef}. In \cite{TT}, the NLFT of $f$ or $\{c_n\}$ is defined as the limit of $ \left( a(N, z), b(N, z) \right)$ as $N \to \infty$, if such a limit exists in some sense.  If the discrete measure $f$ is supported on finitely many points, $f = \sum_{n = 1}^N c_n \delta_n$, then the NLFT becomes a finite product of matrix exponentials \begin{equation*}
    \begin{pmatrix} a^\#(\infty, z) & b^\#(\infty, z) \\ b(\infty, z) &  a(\infty, z) \end{pmatrix} = \prod_{n = 1}^N \exp \begin{pmatrix} 0 & c_n \exp( - 2 i z n) \\ c_n \exp(2 i z n) & 0 \end{pmatrix}.
\end{equation*}

Even in the finite product case above, this formula is still only symbolic, as matrix exponentials are hard to compute. However, this is the only example of the NLFT (to our knowledge) that exists in the literature.

For an infinite one-sided sequence (the potential of the real Dirac system is an evenly spaced discrete measure on $\mathbb{R}_+$), it is known that the limit of  $ \left( a(n, z), b(n, z) \right)$ exists when $\{c_n\} \in \ell^p$, $1 \leqslant p \leqslant 2$, see \cite{TT}.

Let an evenly spaced discrete measure $f$ be the potential of the Dirac system, and $\mu$ be the spectral measure of the system with the Neumann boundary condition at $t = 0$, we can use definition \ref{NLFTDef} for the discrete case as well. The function $ \nlhat{f} = b(\infty, z)/ a(\infty, z)$ also satisfies \begin{equation*}
    \nlhat{f} = \frac{b(\infty, z)}{a(\infty, z)} = \frac{\mathcal{S}_\mu - 1}{\mathcal{S}_\mu + 1}.
\end{equation*}

\section{Inverse Spectral Problems}\label{PZ}

The goal of inverse spectral problem is to recover the original differential equation system from its spectral data. This section summarizes some of the recent developments in inverse spectral problems for canonical systems. 

Let $a_n$ and $b_n$ be trigonometric moments of a given locally finite $2\pi$-periodic measure $\mu$,
$$a_0=\frac {\mu([-\pi,\pi])}{2\pi},\ \ a_n=\frac 1{\pi}\int_{-\pi}^{\pi}\cos (nx) \ d\mu(x),\ \ b_n=\frac 1{\pi}\int_{-\pi}^{\pi}\sin (nx) \ d\mu(x),\ n=1,2,...$$
we will formally write $d\mu(x)=\sum_{n = 0}^\infty (a_n \cos(nx) + b_n \sin (nx)) dx$.

The real Dirac systems we use to study the non-linear Fourier transform always give rise to diagonal determinant-normalized canonical systems without jump intervals. Moreover, canonical systems with diagonal Hamiltonians always give rise to even spectral measures (see \cite{MP} and references therein for details). The following two results provide an explicit way to recover the Hamiltonian matrix in \eqref{CS} if the spectral measure is even and periodic. Since we only consider systems with even spectral measures, we have $b_n = 0$ for all $n$. 

\begin{theorem} \label{2piThm} \cite{MP}
Let $\mu$ be a spectral measure of a system \eqref{CS}. Suppose that $\mu$ is an even $2\pi$-periodic PW-measure, $d\mu(x) = \sum_{n = 0}^\infty a_n \cos(nx) dx$.

Consider the infinite Toeplitz matrix \begin{equation*}
    J = \begin{pmatrix} a_0 & \frac{a_1}{2} & \frac{a_2}{2} & \frac{a_3}{2}  & \ldots\\
    \frac{a_1}{2} & a_0 & \frac{a_1}{2} & \frac{a_2}{2}  & \ldots\\
    \frac{a_2}{2} & \frac{a_1}{2} & a_0 & \frac{a_1}{2}  & \ldots\\
    \frac{a_3}{2} & \frac{a_2}{2} & \frac{a_1}{2} & a_0  & \ldots\\
    \vdots & \vdots & \vdots & \vdots   & \ddots\\
    \end{pmatrix}.
\end{equation*}

Denote by $J_n$ the $(n + 1) \times (n + 1)$ sub-matrix on the upper-left corner of $J$. Then $h_{11}(t)$ in the Hamiltonian is a step function with \begin{equation*}
    h_{11}(t) = \Sigma(J_n^{-1}) - \Sigma(J_{n - 1}^{-1}) \; \text{on} \; \Big(\frac{n}{2}, \frac{n + 1}{2}\Big],\ n = 0, 1, 2, \ldots
\end{equation*}
where $\Sigma$ denotes the sum of all elements in the matrix, and $J_{- 1}^{-1} = 0$.
\end{theorem}

Via a change of variable one can extend the last statement to $2T-$periodic measures:

\begin{corollary}\label{2Talgorithm}\cite{PZ}
Let $\mu$ be an even $2T$-periodic PW-measure, $d\mu(x) = \left( \sum_{n = 0}^\infty a_n \cos\left(\frac{2 n \pi}{2T} x\right) \right)dx$. Consider the infinite Toeplitz matrix \begin{equation*}
    J = \begin{pmatrix} a_0 & \frac{a_1}{2} & \frac{a_2}{2} & \frac{a_3}{2} & \ldots\\
    \frac{a_1}{2} & a_0 & \frac{a_1}{2} & \frac{a_2}{2} & \ldots\\
    \frac{a_2}{2} & \frac{a_1}{2} & a_0 & \frac{a_1}{2}  & \ldots\\
    \frac{a_3}{2} & \frac{a_2}{2} & \frac{a_1}{2} & a_0 & \ldots\\
    \vdots & \vdots & \vdots & \vdots & \ddots\\
    \end{pmatrix}.
\end{equation*}

Denote by $J_n$ the $(n + 1) \times (n + 1)$ sub-matrix on the upper-left corner of $J$. Then $h_{11}(t)$ in the Hamiltonian is a step function with \begin{equation*}
    h_{11}(t) = \Sigma(J_n^{-1}) - \Sigma(J_{n - 1}^{-1}) \; \text{on} \; \Big(\frac{n\pi}{2T}, \frac{(n + 1)\pi}{2T}\Big], n = 0, 1, 2, \ldots
\end{equation*}
where $\Sigma$ denotes the sum of all elements in the matrix, and $J_{-1}^{-1} = 0$.
\end{corollary}

Theorem \ref{2piThm} can be rewritten using the values of orthonormal polynomials on the unit circle, thereby establishing a connection between spectral problems and the theory of orthogonal polynomials on the unit circle (OPUC).

\begin{theorem}\label{onpoly}\cite{MP}
The value of the $n$-th step of $h_{11}(t)$ is $\abs{\varphi_n(1)}^2$, where $\varphi_n(z)$ is the $n$-th orthonormal polynomial on $\mathbb{T}$ with $\mu_\mathbb{T} = \frac{1}{2\pi} \mu|_{[-\pi, \pi)}$.
\end{theorem}

For non-periodic spectral measures, the inverse spectral problem can be solved using a "periodization" approach. Let $\mu$ be an even spectral measure of a canonical system. We denote by $\mu_T$ its $2T-$periodization, defined as $\mu_T(S) = \mu(S)$ for $S\subseteq [-T, T)$, and extended periodically to the rest of $\mathbb{R}$.

\begin{theorem}\cite{PZ}\label{PWConverge}
Let $\mu$ be an even spectral measure of a PW-system, and $\mu_T$ be its periodizations, then $h_{11}^T \overset{\ast}{\to} h_{11}$ as $T \to \infty$.
\end{theorem}

Here, $h_{11}^T \overset{\ast}{\to} h_{11}$ means $\int_0^\infty h_{11}^T(t)\phi(t)dt{\to} \int_0^\infty h_{11}(t)\phi(t)dt$ as $T \to \infty$, for every continuous compactly supported function $\phi$ on $\mathbb{R}_+$.

\section{Main Results}\label{results}

If a discrete measure supported on positive half-integers is used as the potential of a real Dirac system, the corresponding canonical system has a diagonal Hamiltonian, where $h_{11}$ is a step function with a step size of $1/2$. This is the case considered in theorems \ref{2piThm} and \ref{onpoly}. We identify a $2\pi-$periodic measure $\mu$ on the real line $\mathbb{R}$ with the measure $\mu = \frac{1}{2\pi} \mu|_{[-\pi, \pi)}$ on the unit circle $\mathbb{T}$. Via direct computation, theorem \ref{onpoly} yields the following result:

\begin{theorem}\label{discreteINLFT}
Let $\mu_\pi$ be an even $2\pi-$periodic PW-measure, and let $\{\varphi_n\}_{n = 0}^\infty$ be the family of polynomials orthonormal in $L^2_{\mu_\pi}(\mathbb{T})$. The inverse non-linear Fourier transform of $\frac{S_{\mu_\pi} - 1}{S_{\mu_\pi} + 1}$ is given by \begin{equation*}
    \sum_{n = 1}^\infty \log \left( \frac{|\varphi_{n}(1)|}{|\varphi_{n - 1}(1)|} \right) \delta_{n/2}.
\end{equation*}
\end{theorem}

Let the discrete measure $f = \sum_{n = 1}^\infty c_n \delta_{n/2}$ be the potential of a real Dirac system, and $\mu_\pi$ be the $2\pi-$periodic spectral measure of the system. Denote by $w(x)$ the density of the absolutely continuous part of the measure $\mu$. If $\{c_n\} \in \ell_2$, it is known that $w(x)$ satisfies the Szeg\"o condition: $\log \abs{w(x)} \in L^1(\Pi)$, see \cite{Denisov} for details on this and related results.

If an even $2\pi-$periodic PW-measure $\mu_\pi$ satisfies the Szeg\"o condition, we can take this measure to be the spectral measure of a real Dirac system, take the inverse NLFT of $\frac{\mathcal{S}_\mu - 1}{\mathcal{S}_\mu + 1}$ using the theorem above, and obtain a discrete measure supported on positive half-integers. It can be shown that this sequence of masses needs to be square summable, which is stated in the following corollary (see \cite{TTT} for details).

\begin{corollary}
Let $\mu_\pi$ be an even $2\pi-$periodic PW-measure satisfying the Szeg\"o condition, and let $\{\varphi_n\}_{n = 0}^\infty$ be the family of polynomials orthonormal in $L^2_{\mu_\pi}(\mathbb{T})$. The sequence \begin{equation*}
    \log \left( \frac{|\varphi_{n}(1)|}{|\varphi_{n - 1}(1)|} \right)
\end{equation*}
is in $\ell_2$.
\end{corollary}

It was mentioned in \cite{TT} that one can use the NLFTs of sequences to approximate the NLFTs of functions and obtain a convergence result. Our next goal is to use this idea to get a convergence result for the inverse NLFT.

Let $\mu$ be the spectral measure of a real Dirac system, and $\mu_T$ be the $2T-$periodization of $\mu$, defined by $\mu_T(S) = \mu(S)$ if $S \subseteq [-T, T)$ and $\mu_T$ is $2T-$periodic. We study the inverse NLFT of the functions $$f(t) = \frac{\mathcal{S}_\mu - 1}{\mathcal{S}_\mu + 1}$$ and $$f^T(t) = \frac{\mathcal{S}_{\mu_T} - 1}{\mathcal{S}_{\mu_T} + 1}.$$ Using results from section \ref{PZ}, we can solve the inverse spectral problems and obtain determinant-normalized diagonal Hamiltonians from $\mu$ and $\mu_T$'s. Denote by $h_{11}$ and $h_{11}^T$ the first entry in the Hamiltonian corresponding to $\mu$ and $\mu_T$, respectively. Then we have \begin{equation*}
\begin{split}
    \nlihat{f}(t) &= \frac{1}{2} \frac{d}{dt} \log \left( h_{11}(t) \right), \\
    \nlihat{f^T}(t) &= \frac{1}{2} \frac{d}{dt} \log \left( h_{11}^T(t) \right),    
\end{split}    
\end{equation*}
where $\nlihat{f}$ and $\nlihat{f^T}$ stand for the inverse NLFT of $f$ and $f^T$, respectively.

Functions from the closed unit ball in $H^\infty(\HP)$ (or $H^\infty(\D)$) are called \textit{Schur functions}. Based on definition \ref{NLFTDef}, the NLFTs of evenly spaced discrete measures on the half-line $\mathbb{R}_+$ and $L^p$ functions on the half-line $\mathbb{R}_+$ are Schur functions.

\begin{theorem}
Let $f(z)$ be a Schur function corresponding to the spectral measure $\mu$ of a Dirac system. Recall that\begin{equation*}
    f(z) = \frac{\mathcal{S}_\mu(z) - 1}{\mathcal{S}_\mu(z) + 1}.
\end{equation*}

Denote by $f^T(z)$, the $2T-$periodization of $f(z)$, defined as the following: Let $\mu_T$ be the $2T-$periodization of $\mu$, and \begin{equation*}
    f^T(z) = \frac{\mathcal{S}_{\mu_T}(z) - 1}{\mathcal{S}_{\mu_T}(z) + 1}.
\end{equation*}

Then $f^T(z) \to f(z)$ pointwise on $\mathbb{R}$, and normally on $\mathbb{C}_+$.

Denote by $\nlihat{f}$ and $\nlihat{f^T}$ the inverse NLFT of $f$ and $f^T$, respectively. Then we also have \begin{equation}\label{weakconv}
    \exp( 2 \int_0^t \nlihat{f^T}(x) dx) \overset{\ast}{\to} \exp( 2 \int_0^t \nlihat{f}(x) dx).
\end{equation}
\end{theorem}

Note that convergence of $\nlihat{f^T}$ to $\nlihat{f}$ follows from theorem \ref{PWConverge}. It remains to prove $f^T \to f$ on $\mathbb{R}$ and normally on $\mathbb{C}_+$.

\begin{proof}

It follows from the definition of $\mu_T$ and theorem \ref{PWcondition} that for a PW-measure $\mu$ and a constant $C$ such that $\mu(x, x + 1) < C$ for all $x \in \mathbb{R}$, we also have that $\mu_T(x, x + 1) < 2C$ for all $x \in \mathbb{R}$ and $T > 0$.

We first prove pointwise convergence on $\mathbb{R}$: fix $z \in \mathbb{R}$: \begin{equation*}
\begin{split}
    \abs{S_\mu(z) - S_{\mu_T}(z)} &= \abs{\frac{1}{\pi i} \int_\mathbb{R} \frac{1}{t - z} - \frac{t}{1 + t^2} d\mu(t) - \frac{1}{\pi i} \int_\mathbb{R} \frac{1}{t - z} - \frac{t}{1 + t^2} d{\mu_T}(t)} \\
    &= \abs{\frac{1}{\pi i} \int_{ \abs{t} > T} \frac{1}{t - z} - \frac{t}{1 + t^2} d(\mu - \mu_T)(t)}\\
    &\leqslant \abs{\frac{1}{\pi i}} \int_{ \abs{t} > T} \abs{\frac{1}{t - z} - \frac{t}{1 + t^2} }d\abs{\mu - \mu_T}(t)\\   
    &\leqslant \frac{6C}{\pi} \int_{ \abs{t} > T} \abs{\frac{1}{t - z} - \frac{t}{1 + t^2} }dt \to 0, \quad \text{as} \; T \to \infty.
\end{split}
\end{equation*}

Let $K$ be a compact subset in $\mathbb{C}_+$, then there exists $R > 0$ such that $|z| < R$ for all $z \in K$. Note that $S_\mu(z) = P_\mu(z) + i Q_\mu(z)$, where $\mathcal{P}_\mu(z)$ is the Poisson integral \begin{equation*}
    \mathcal{P}_\mu(z) = \mathcal{P}_\mu(x + iy) = \frac{1}{\pi} \int_\mathbb{R} \frac{y}{(t - x)^2 + y^2} d\mu(t),
\end{equation*} and $\mathcal{Q}_\mu(z)$ is the conjugate Poisson integral \begin{equation*}
    \mathcal{Q}_\mu(z) = \mathcal{Q}_\mu(x + iy) = \frac{1}{\pi} \int_\mathbb{R} \frac{x - t}{(x - t)^2 + y^2} + \frac{t}{1 + t^2} d\mu(t).
\end{equation*} 
It suffices to show uniform convergence of the Poisson and conjugate Poisson integrals on $K$. We start with the Poisson integrals. Let $T > R$, \begin{equation*}
\begin{split}
    \abs{P_\mu(z) - P_{\mu_T}(z)} & = \abs{\frac{1}{\pi} \int_{\mathbb{R}} \frac{y}{(t - x)^2 + y^2} d\mu(t) - \frac{1}{\pi} \int_{\mathbb{R}} \frac{y}{(t - x)^2 + y^2} d\mu_T(t)}\\
    &= \frac{1}{\pi} \abs{\int_{\abs{t} > T} \frac{y}{(t - x)^2 + y^2} d (\mu - \mu_T)(t) }\\
    &\leqslant \frac{6C}{\pi} \int_{\abs{t} > T} \abs{\frac{y}{(t - x)^2 + y^2}} dt \leqslant \frac{6CR}{\pi} \int_{\abs{t} > T} \frac{1}{(t - x)^2} dt \\
    &\leqslant \frac{6CR}{\pi} \int_{\abs{t} > T} \frac{1}{(t - R)^2}dt \to 0, \quad \text{as} \; T \to \infty.
\end{split}
\end{equation*}

Finally, we look at the conjugate Poisson integrals.
\begin{equation*}
\begin{split}
    \abs{Q_\mu(z) - Q_{\mu_T}(z)} & = \abs{\frac{1}{\pi} \int_{\mathbb{R}} \left( \frac{x - t}{(x - t)^2 + y^2} + \frac{t}{1 + t^2} \right) d(\mu - \mu_T)(t)}\\
    &= \frac{1}{\pi} \abs{\int_{\abs{t} > T} \left( \frac{x - t}{(x - t)^2 + y^2} + \frac{t}{1 + t^2} \right) d(\mu - \mu_T)(t)}\\
    &\leqslant \frac{6C}{\pi} \int_{\abs{t} > T} \abs{ \frac{x - t}{(x - t)^2 + y^2} + \frac{t}{1 + t^2}} dt.
\end{split}
\end{equation*}
For every fixed $z = x + iy \in K$, $\int_{\abs{t} > T} \abs{ \frac{x - t}{(x - t)^2 + y^2} + \frac{t}{1 + t^2}} dt$ is monotonically decreasing in $T$, and converges to $0$. By Dini's theorem, $\int_{\abs{t} > T} \abs{ \frac{x - t}{(x - t)^2 + y^2} + \frac{t}{1 + t^2}} dt$ converges uniformly on $K$. 
\end{proof}

\begin{remark}
Convergence established in the theorem above from $\nlihat{f^T}$ and $\nlihat{f}$ is weak. We are able to obtain a stronger convergence result for certain examples. See example \ref{strongconv} below. However, we cannot strengthen the mode of convergence in general because of the following examples.
\end{remark}

\begin{example}
This example shows \eqref{weakconv} is not preserved in general if logarithm is taken on both sides.

Let $f_0(t)$ be the $1$-periodic step function on $\R_+$ with \begin{equation*}
    f_0(t) = \begin{cases} 1 \quad &t \in \left[ 0, \frac{1}{2} \right)\\
    \frac{1}{4}, &t \in \left[\frac{1}{2}, 1 \right)
    \end{cases}.
\end{equation*}

Let $f^T(t) = f_0(T \cdot t)$, then $f^T(t)$'s converge weakly to the constant function $f = \left( 1 + \frac{1}{4} \right)/2 = \frac{5}{8}$. If we take the logarithm on both $f^T$ and $f$, we get $$\log \left( f_T(t) \right) \overset{\ast}{\to} \left( \log(1) + \log \left( \frac{1}{4} \right) \right)/2 \neq \log\left( f\right).$$

We cannot eliminate examples like this. These $f^T$'s, if taken to be $h_{11}(t)$'s in canonical systems \eqref{CS} with diagonal determinant-normalized Hamiltonians, give rise to periodic spectral measures $\mu_T$'s. Whether these $\mu_T$'s converge to $\mu$, the spectral measure of the canonical system with $h_{11}(t) = f$, is still unclear. If $\mu_T \to \mu$ in some sense, this example would suggest that the limit of discrete NLFTs might not be the NLFT of the limit of the discrete measures, implying that \eqref{weakconv} cannot be strengthened.
\end{example}

\begin{remark}  
If we think of the $2$ inside the exponential in \eqref{weakconv} as a square, then we can construct a similar example to show that we cannot take the squareroot while preserving weak convergence. 
\end{remark}

\section{Examples of the Non-linear Fourier Transform}\label{examples}

The following examples of the NLFT are derived from solutions to inverse spectral problems for canonical systems in the work of Makarov and Poltoratski \cite{MP}.

\begin{example}
    \begin{equation*}
    f(t) = \sum_{k = 1}^\infty \log \left( \frac{k + 2}{k} \right) \delta_{k/2} ;
\end{equation*}

\begin{equation*}
    \nlhat{f}(z) = \frac{- e^{iz}}{2 - e^{iz}} ;
\end{equation*}

\begin{equation*}
    a(\infty, z) = \frac{2 - e^{iz}}{\sqrt{2}(1 - e^{iz})}, \quad b(\infty, z) = \frac{- e^{iz}}{\sqrt{2}(1 - e^{iz})}.
\end{equation*}
\end{example}

\begin{example}
    \begin{equation*}
    f(t) = C \delta_1;
\end{equation*}

\begin{equation*}
    \nlhat{f}(z) = \tanh{(2C)} \cdot e^{iz/2} ;
\end{equation*}

\begin{equation*}
    a(\infty, z) = \cosh{(2C)}, \quad b(\infty, z) = \sinh{(2C)} \cdot e^{iz/2}.
\end{equation*}
    
\end{example}

\begin{remark}
This example shows that the NLFT behaves similarly to its linear counterpart: a point mass is transformed to an exponential. However, the constant factor in front is transformed non-linearly.
\end{remark}

\begin{example}\label{strongconv}
Consider the spectral measure $\mu = \sqrt{2\pi} \delta_0 + \frac{1}{\sqrt{2\pi}}m$, where $\delta_0$ is the unit point mass at the origin, and $m$ is the Lebesgue measure.

We first compute $f$ and $\nlhat{f}$ from $\mu$. As shown in \cite{MP}, $h_{11}(t) = \frac{\sqrt{2\pi}}{(1 + 2t)^2}$, then $f(t) = \frac{-4 \sqrt{2 \pi}}{1 + 2t}$. We also have \begin{equation*}
    \nlhat{f}(z) = \frac{\mathcal{S}_\mu - 1}{\mathcal{S}_\mu + 1} = \frac{-2 \sqrt\pi z + \sqrt2 z + 2i \sqrt2}{2 \sqrt\pi z + \sqrt2 z + 2i \sqrt2} = \frac{\left(1 - 2\pi\right)z^2 + 4i\sqrt{2\pi}z + 4}{\left(1 + 2\pi + 2 \sqrt{2\pi}\right)z^2 + 4}.
\end{equation*}

The $2T$-periodization $\mu_T = \frac{1}{\sqrt{2\pi}} m + \sum_{n \in \mathbb{Z}} \sqrt{2\pi} \delta_{2 n T}$. As shown in \cite{PZ}, the $n-$th step of $h_{11}^T$ takes value \begin{equation*}
    \frac{\sqrt{2\pi} T^2}{(n \pi + T)(n \pi + T + \pi)}.
\end{equation*}
Then by theorem \ref{discreteINLFT} the corresponding potentials in the Dirac systems are  \begin{equation*}
    f(t) = \sum_{n \in \mathbb{Z}_+} \log \left( \frac{n \pi + T}{n \pi + 2\pi + T} \right) \delta_{\frac{n\pi}{2T}},
\end{equation*}
the NLFTs of these $f^T$'s are \begin{equation*}
    \nlhat{f^T}(z) = \frac{\mathcal{S}_{\mu_T} - 1}{\mathcal{S}_{\mu_T} + 1} = \frac{i \left( \sqrt{2} \pi \cot (\pi z /2T) + 2i \sqrt{\pi} T - i \sqrt{2} T \right)}{i \sqrt{2} \pi \cot(\pi z /2T) + 2 \sqrt{\pi} T + \sqrt{2} T}.
\end{equation*}

Now let $T \to \infty$, we have \begin{equation*}
    \lim_{T \to \infty} \nlhat{f^T}(z) =  \frac{\left(1 - 2\pi\right)z^2 + 4i\sqrt{2\pi}z + 4}{\left(1 + 2\pi + 2 \sqrt{2\pi}\right)z^2 + 4},
\end{equation*}

For this example, we get convergence stronger than \eqref{weakconv}. If the discrete measures $f^T = \sum_{n \in \mathbb{Z}_+} a_n^T \delta_{\frac{n\pi}{2T}} \to f$ $\ast-$weakly, we expect the values $a_n^T/\frac{\pi}{2T}$ getting closer to $f(\frac{n\pi}{2T})$, which is what we see in figure \ref{ftf}. 

\begin{figure}
\centering
\captionsetup{justification=centering}
\includegraphics[width=0.45\textwidth]{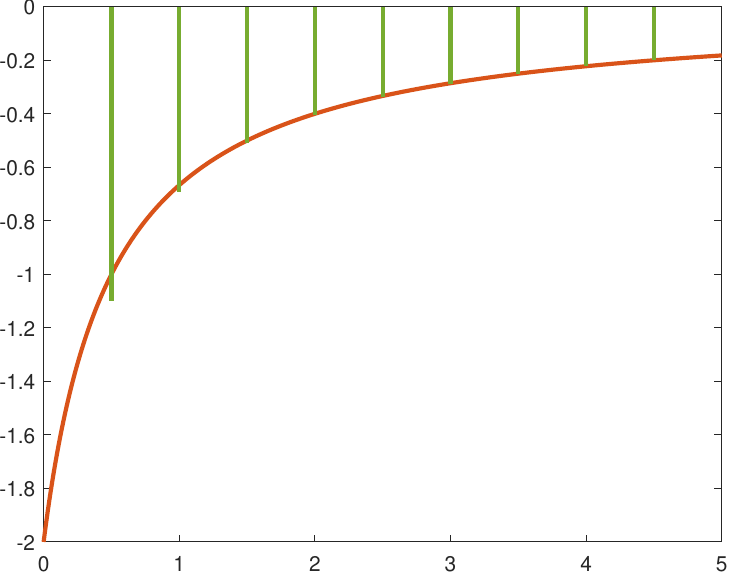}
\includegraphics[width=0.45\textwidth]{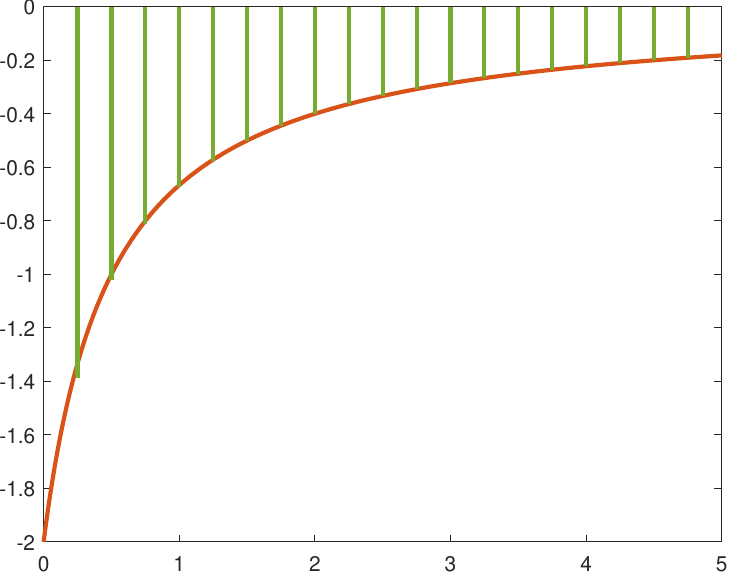}
\includegraphics[width=0.45\textwidth]{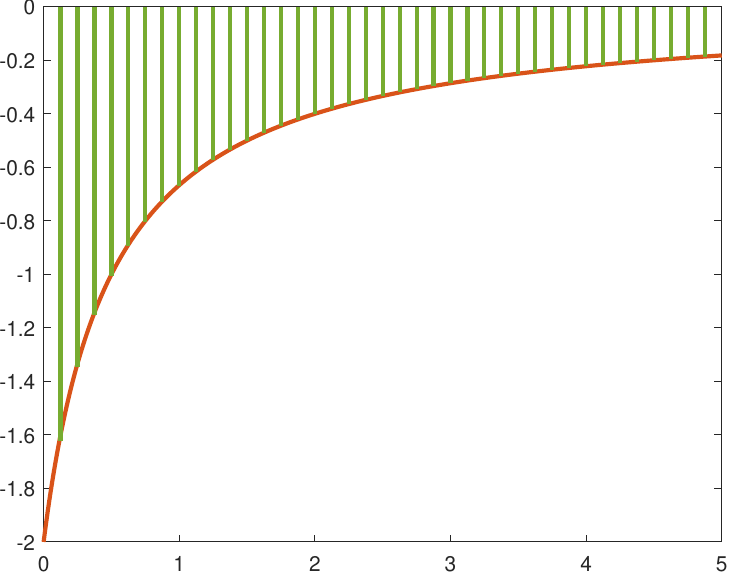}
\includegraphics[width=0.45\textwidth]{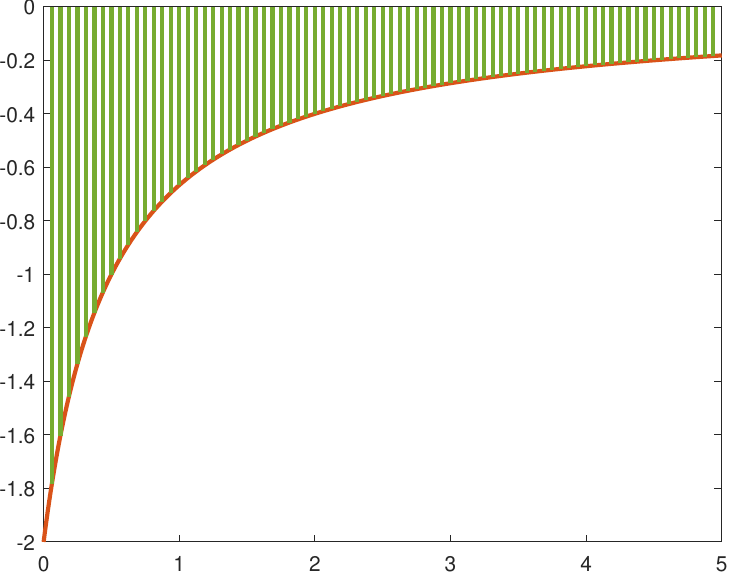}
\caption{Convergence of the inverse NLFT}
\caption*{First row from left to right: $T = \pi$, $T = 2\pi$.\\ Second row from left to right: $T = 4\pi$, $T = 8\pi$.\\ Orange continuous curves: the function $\nlihat{f}$.\\ Green line segments: $a_n^T \frac{2T}{\pi}$ at $\frac{n\pi}{2T}$.  }
\label{ftf}
\end{figure}

\end{example}

\begin{example}
We can further generalize the measure in the example above to the following: let $\mu = m + \beta \pi \delta_0$, then 
    \begin{equation*}
    f(t) = - \frac{\beta}{1 + \beta t}
\end{equation*}
and
\begin{equation*}
    \nlhat{f}(z) = \frac{\beta i}{2z + \beta i}.
\end{equation*}

The traditionally studied $(a, b)-$pair can be recovered from $\nlhat{f}$ as well:
\begin{equation*}
    a(\infty, z) = \frac{2z + \beta i}{2z}, \quad b(\infty, z) = \frac{\beta i}{2z}.
\end{equation*}
\end{example}

\newpage
\printbibliography

\end{document}